\documentclass[a4paper, 12pt, oneside]{article} 
\usepackage{graphicx} 
\usepackage{amsmath} 
\usepackage{amssymb} 
\usepackage{amsfonts} 
\usepackage{amsthm} 
\usepackage{longtable} 
\usepackage{multirow}
\usepackage{color}
\usepackage{makeidx}

\makeindex

\title{Icosahedral Skeletal Polyhedra Realizing Petrie Relatives of Gordan's Regular Map}
\author{Anthony M. Cutler\\
Boston, MA 02115, USA\\[.08in]
{\small and} \\[.05in]
Egon Schulte\thanks{Supported by NSF-grant DMS--0856675}\\
Northeastern University\\
Boston, MA 02115, USA\\[.05in]
{\small and} \\[.08in]
J\"org M. Wills\\
University of Siegen\\
D-57068 Siegen, Germany}

\date{\sl\small\today }

\newtheorem{lemma}{Lemma}[section] 
\newtheorem{theorem}[lemma]{Theorem}

\renewcommand{\Gamma}{\varGamma} 
\renewcommand{\epsilon}{\varepsilon}

\newcommand{\E}{\mathbb{E}^3}

\begin{document}

\maketitle

\begin{abstract}
\noindent
Every regular map on a closed surface gives rise to generally six regular maps, its {\em Petrie relatives\/}, that are obtained through iteration of the duality and Petrie operations (taking duals and Petrie-duals). It is shown that the skeletal polyhedra in Euclidean $3$-space which realize a Petrie relative of the classical Gordan regular map and have full icosahedral symmetry, comprise precisely four infinite families of polyhedra, as well as four individual polyhedra.
\end{abstract}

\section{Introduction}
\label{intro}

It is well-known that the Petrie relatives of a regular map on a closed surface form a family of generally six regular maps with the same automorphism group, obtained from the original map under iteration of the duality and Petrie operations (taking duals and Petrie-duals); see Coxeter \& Moser~\cite{cm}, McMullen \& Schulte~\cite[Ch. 7B]{arp} or Wilson~\cite{wilson}. The purpose of this short note is to point out that the skeletal polyhedra in Euclidean $3$-space $\E$ which realize a Petrie relative of the classical Gordan~\cite{gor} regular map and have full icosahedral symmetry, comprise precisely four infinite families of polyhedra, as well as four additional, individual polyhedra. Here a skeletal polyhedron is a finite geometric edge-graph in space equipped with a polyhedral face structure; see Gr\"unbaum~\cite{gr1} and \cite[Ch. 7E]{arp}, as well as Section~\ref{maskpo}.

The icosahedral polyhedra described in this paper are particularly interesting examples of regular polyhedra of index $2$ (see Cutler~\cite{cut}, Cutler \& Schulte~\cite{cutsch}, Wills~\cite{wills}); that is, they are combinatorially regular but ``fail geometric regularity by a factor of $2$". Four polyhedra were previously known and have an interesting history; namely, the planar-faced polyhedral realizations of Gordan's map $\{5,4\}_6$ and its dual $\{4,5\}_6$ of genus $4$, as well as of $\{6,5\}_4$ and $\{5,6\}_4$ of genus $9$. These give four of just five regular polyhedra of index $2$ which are orientable and have planar faces (see \cite{wills} for the enumeration and Richter~\cite{rich} for figures). The polyhedra for $\{5,4\}_6$ and $\{5,6\}_4$ were discovered by Hess (1878), Pitsch (1881) and Badoureau (1881); for the history see Coxeter~\cite[\S  6.4]{coxeter}, and for figures see \cite[Fig. 4.4a]{coxeter} and Coxeter, Longuet-Higgins \& Miller~\cite[Figs. 45, 53]{clm}. The polyhedra for $\{4,5\}_6$ and $\{6,5\}_4$, respectively, are obtained from the polyhedra for $\{5,4\}_6$ and $\{5,6\}_4$ by polarity (see Gr\"unbaum \& Shephard~\cite{gsh}); for figures see \cite[Fig. 6.4c]{coxeter}, \cite[$De_{2}f_{2}$ on Plate XI]{cdfp}, and \cite{swirp}. 

In addition to these classical examples, new icosahedral polyhedra for Petrie relatives of Gordan's map can be derived from the complete classification of the finite regular polyhedra of index $2$ in \cite{cut,cutsch}. It turns out that this provides a complete list of possible realizations. Overall, each Petrie relative is seen to admit a skeletal realization with icosahedral symmetry, and some relatives can even be realized in more than one way.

\section{Maps and skeletal polyhedra}
\label{maskpo}

Skeletal polyhedra were first investigated in Gr\"unbaum~\cite{gr1}. A priori they are not solid figures bounded by faces spanned by membranes, but rather more general discrete polyhedral structures with convex or non-convex, planar or skew, and finite or infinite (helical or zig-zag) polygonal faces and vertex-figures. Here we are only concerned with finite polyhedral structures, so in particular helical or zigzag faces will not occur. 

To begin with, a (finite) {\em polygon\/}, or simply an {\em $n$-gon\/}, $(v_1, v_2, \dots, v_n)$ in $\E$ is a figure formed by distinct points $v_1, \dots, v_n$, together with the line segments $(v_i, v_{i+1})$, for $i = 1, \dots, n-1$, and $(v_n, v_1)$. The points and line segments are the {\em vertices\/} and {\em edges\/} of the polygon. Note that polygons are edge cycles, {\em not\/} topological discs, and may or may not be plane polygons. 

A (finite) {\em skeletal polyhedron\/}, or simply {\em polyhedron\/}, $P$ in $\E$ consists of a finite set $V_P$ of points, called {\em vertices\/}, a finite set $E_P$ of line segments, called {\em edges\/}, joining points of $V_P$, and a finite set $F_P$ of polygons, called {\em faces\/}, formed by line segments of $E_P$ such that the following three conditions hold. First, $P$ is globally connected, in the sense that the {\em edge graph\/} of $P$ defined by $V_P$ and $E_P$ is connected; and second, $P$ is locally connected, meaning that the vertex-figure of $P$ at every vertex of $P$ is connected. Recall here that the {\em vertex-figure\/} of $P$ at a vertex $v$ is the graph whose vertices are the neighbors of $v$ in the edge graph of $P$ and whose edges are the line segments $(u,w)$, where $(u, v)$ and $(v, w)$ are edges of a common face of $P$. Third, each edge of $P$ is an edge of exactly two faces of~$P$. 

For combinatorial purposes we identify a polyhedron $P$ with the underlying map on a closed compact (orientable or non-orientable) surface; see \cite{cm,arp}. Thus skeletal polyhedra in $\E$ are $3$-dimensional geometric realizations of abstract polyhedra in the sense of \cite[Ch. 7E]{arp} (see also \cite{gr1,grhol} and \cite{ordinary}). Note, however, that a skeletal polyhedron is a finite geometric edge graph in $\E$ equipped with a distinguished face structure; it is not a solid figure.  

A map $P$ is ({\em combinatorially\/}) {\em regular\/} if its {\em automorphism group\/} $\Gamma(P)$ is transitive on the flags (incident vertex-edge-face triples). The automorphism group $\Gamma(P)$ of a regular map is generated by involutions $\rho_0,\rho_1,\rho_2$ satisfying (at least) the Coxeter-type relations
\begin{equation}
\label{relone}
\rho_{0}^{2} = \rho_{1}^{2} = \rho_{2}^{2} =
(\rho_{0}\rho_{1})^{p} = (\rho_{1}\rho_{2})^{q} = (\rho_{0}\rho_{2})^{2} = 1,
\end{equation}
where $p$ and $q$ determine the {\em type} $\{p,q\}$ of $P$. Thus $\Gamma(P)$ is a quotient of the Coxeter group $[p,q]$ abstractly defined by the relations in (\ref{relone}). 

The {\em Petrie-dual\/} of a regular map $P$ is a (regular) map with the same vertices and edges as $P$, obtained by replacing the faces of $P$ by the Petrie polygons of $P$. Recall that a {\em Petrie polygon\/} of $P$ is a polygon along the edges of $P$ such that any two, but no three, consecutive edges belong to a common face. The iteration of taking duals and Petrie-duals gives rise to a family of generally six regular maps, the {\em Petrie relatives\/} of $P$, all sharing the same automorphism group. 

Following \cite{cm}, we let $\{p,q\}_r$ denote a regular map derived from the spherical, Euclidean, or hyperbolic regular plane tessellation $\{p,q\}$ by identifying any two vertices~$r$ steps apart along a Petrie polygon of $\{p,q\}$. The  automorphism group of $\{p,q\}_r$ is abstractly defined by (\ref{relone}) and the relation $(\rho_{0}\rho_{1}\rho_{2})^{r}=1$. The dual and Petrie-dual of $\{p,q\}_r$ are the maps $\{q,p\}_r$ and $\{r,q\}_p$, respectively. Thus the Petrie relatives of $\{p,q\}_r$ are given by $\{p,q\}_{r}$, $\{q,p\}_{r}$, $\{r,q\}_{p}$, $\{q,r\}_{p}$, $\{r,p\}_{q}$, and $\{p,r\}_{q}$. Note that $\{p,q\}_r$ is orientable if and only if $r$ is even. 

The present paper focuses on icosahedral polyhedra realizing a Petrie relative of Gordan's map $\{5,4\}_6$. In this case the six Petrie relatives are given by 
\begin{equation}
\label{gord}
\{5,4\}_{6},\, \{4,5\}_{6},\, \{6,4\}_{5},\ \{4,6\}_{5},\, \{6,5\}_{4},\, \{5,6\}_{4} ,
\end{equation}
and all have the same automorphism group of order $240$. The first two maps are orientable of genus $4$; the next two are non-orientable of genus $12$; and the last two are orientable of genus $9$.

A polyhedron with full icosahedral symmetry realizing a map in (\ref{gord}) must be a regular polyhedron of index~$2$ in the sense of~\cite{wills,cutsch}. Recall that $P$ is a {\em regular polyhedron of index~$2$\/} if $P$ is a regular map and the geometric symmetry group $G(P)$ of $P$ is a subgroup of $\Gamma(P)$ of index $2$. For Gordan's map and its relatives, $G(P)$ is the full icosahedral group of order $120$ and has index $2$ in $\Gamma(P)$.

\section{The polyhedra}
\label{pols}

It follows from the above analysis that the desired icosahedral skeletal polyhedra for Gordan's map and its relatives are precisely the regular polyhedra of index $2$ that have full icosahedral symmetry and are combinatorially isomorphic to a map in (\ref{gord}).
This allows us to draw upon the complete enumeration of the finite regular polyhedra of index~$2$ obtained in \cite{cut,cutsch}.  
It was shown in \cite{cutsch} that up to similarity there are precisely $22$ infinite families of regular polyhedra of index $2$ with vertices on two orbits under the full symmetry group, where two polyhedra belong to the same family if they differ only in the relative size of the spheres containing their vertex orbits; all polyhedra in these $22$ families are orientable, but only two have planar faces. The complete enumeration of the remaining regular polyhedra of index $2$, with vertices on one orbit, is described in \cite{cut}; there are exactly $10$ such (individual) polyhedra. 

In particular, this establishes the following theorem.

\begin{theorem}
\label{theo1}
The skeletal polyhedra in $\E$ which realize a Petrie relative of the Gordan regular map and have full icosahedral symmetry, comprise precisely four infinite families of polyhedra, as well as four individual polyhedra.
\end{theorem}

\begin{table}[htp]
\centering
{\begin{tabular}{|c|c|c|l|c|c|c|cl}  \hline
Type &f-Vector  &Vertex Orbits&Map &Planar Faces\\
$\{p,q\}_r$ &$(f_{0},f_{1},f_{2})$&& of \cite{con}&\\[.05in]
\hline
\hline
$\{4,5\}_{6}$  & $(24,60,30)$ &2&$R4.2$&Yes, for one polyhedron\\
\hline
$\{6,5\}_{4}$  & $(24,60,20)$&2& $R9.16^*$ &Yes, for one polyhedron\\
\hline
$\{4,5\}_{6}$  & $(24,60,30)$&2& $R4.2$& No\\
\hline
$\{6,5\}_{4}$  & $(24,60,20)$&2 & $R9.16^*$& No\\
\hline
\hline
$\{4,6\}_{5}$	&$(20, 60, 30)$& 1&$N12.1$& No\\
\hline
$\{5,6\}_{4}$	&$(20, 60, 24)$&1&$R9.16$&Yes \\
\hline
$\{6,4\}_{5}$	&$(30, 60, 20)$&1&$N12.1^*$&No\\
\hline
$\{5,4\}_{6}$	&$(30, 60, 24)$&1&$R4.2^*$& Yes \\
\hline
\end{tabular}
\caption{The icosahedral skeletal polyhedra realizing Petrie relatives of Gordan's map. The first four rows represent the infinite families of polyhedra shown in Figure~\ref{figcsw1}, the last four the individual polyhedra shown in Figure~\ref{figcsw2}.}
\label{tabone}}
\end{table}

Details about these polyhedra are summarized in Table~\ref{tabone}; the first four rows concern the polyhedra occurring in an infinite family, and the last four rows the individual polyhedra. The $f$-vector $(f_0,f_1,f_2)$ in the second column records the numbers $f_0$, $f_1$ and $f_2$ of vertices, edges and faces of a map. The third column lists the number of vertex orbits, $1$ or $2$, under the full icosahedral symmetry group. The fourth column gives the name of the map in the notation of Conder \cite{con}. Here ``$R$" or ``$N$", respectively, indicates an orientable or non-orientable regular map; the number before the period is the genus, and an asterisk indicates the dual. 

\begin{figure}[hp]
\vspace{-.8in}
\begin{center}
\includegraphics[width=10.2cm, height=13cm]{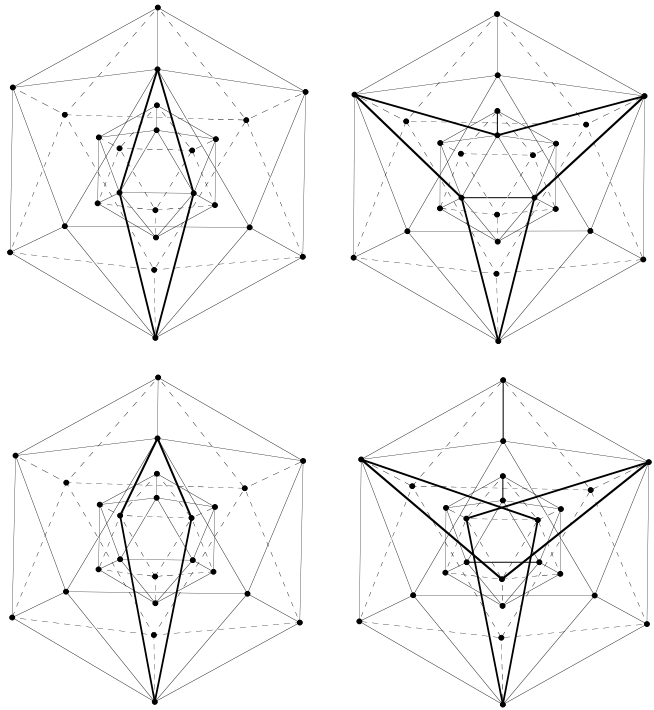}\\[-.55in]
\caption{The four families of icosahedral skeletal polyhedra of types $\{4,5\}_{6}$ or $\{6,5\}_4$.}
\label{figcsw1}
\end{center}
\end{figure}

A member of each of the four families of icosahedral skeletal polyhedra which realize a Petrie relative of the Gordan map is shown in Figure~\ref{figcsw1}. They each are orientable, and have one face orbit, one edge orbit and two vertex orbits. The ratio of the diameters of the vertex orbits, which may be any positive quantity other than 1, identifies the individual polyhedra within each family. The two families in the top row each contain exactly one polyhedron with planar faces, which occurs when the ratio of the diameters of the vertex orbits is $(1 + \sqrt{5})/2$ or $2 + \sqrt{5}$, respectively. The vertices of the polyhedra coincide with the vertices of two concentric regular icosahedra, and these underlying solids, as well as a representative face of each polyhedron, are diagramed. Polyhedra in the same row are Petrie duals, and those in the same column have the same map, of type $\{4,5\}_6$ or $\{6,5\}_4$ respectively, and these maps are dual to the maps of the polyhedra shown in the bottom row of Figure 2. The 30 faces of each polyhedron in the left column are each centered on an edge of the underlying icosahedra, whereas the 20 faces of each polyhedron in the right column are each centered on a face of the underlying icosahedra. \\[-1in]

\begin{figure}[hp]
\begin{center}
\includegraphics[width=12cm, height=15.8cm]{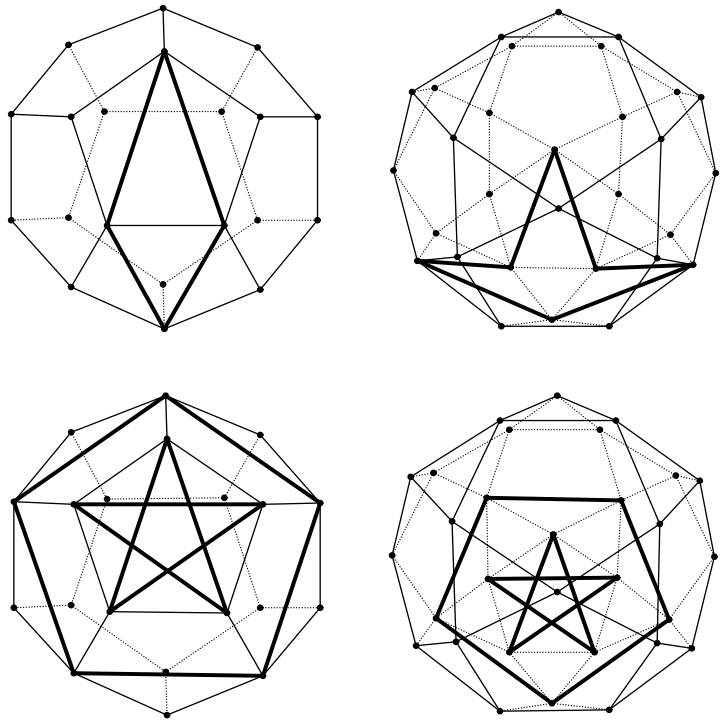}\\[-1.1in]
\caption{The four individual icosahedral skeletal polyhedra which realize a Petrie relative of the Gordan map.}
\label{figcsw2}
\end{center}
\end{figure}

The four individual icosahedral skeletal polyhedra which realize a Petrie relative of the Gordan map are shown in Figure~\ref{figcsw2}. They each have one edge orbit and one vertex orbit. A representative face from each face orbit is diagramed. The vertices coincide with the vertices of either a regular dodecahedron or a regular icosidodecahedron, as shown, and each of the 60 edges of each skeletal polyhedron traverses a pentagonal face of those underlying solids. Polyhedra in the same column are Petrie duals, and the maps of those in the top row are duals. The top left has type $\{4,6\}_5$, is non-orientable and has one face orbit with non-planar faces. The 30 faces are each centered on an edge of the underlying dodecahedron. The bottom left has type $\{5,6\}_4$, is orientable and has two face orbits each with planar faces. The 24 faces are each centered on a face of the underlying dodecahedron, either as an internal pentagram or as a surrounding pentagon. The top right has type $\{6,4\}_5$, is non-orientable and has one face orbit with non-planar faces. The 20 faces are each centered on a triangular face of the underlying icosidodecahedron. The bottom left has type $\{5,4\}_6$, is orientable and has two face orbits each with planar faces. The 24 faces are each centered on a pentagonal face of the underlying icosidodecahedron, again either as an internal pentagram or as a surrounding pentagon.

\end{document}